\documentclass{llncs}

\begin{document}

\title{Modular Forms for $\mbox{GL}(3)$ and Galois Representations}
\author{Bert van Geemen\inst{1} \and Jaap Top\inst{2}}
\institute{Dipartimento di Matematica, Universit\`{a} di Pavia,\\ 
 Via Ferrata 1,  
 27100 Pavia, Italy \\ \email{geemen@dragon.ian.pv.cnr.it}  \and  
IWI, Rijksuniversiteit Groningen, \\ Postbus 800, NL-9700 AV Groningen, The Netherlands
\\ \email{top@math.rug.nl}}
\maketitle
\begin{abstract}
A description and an example are given of numerical experiments which look for a relation
between modular forms for certain congruence subgroups of $\mbox{SL}(3,\bbbz)$
and Galois representations. 
\end{abstract}

\section{Introduction}

In this paper we review a recently  discovered relation between some modular
forms for congruence subgroups of $\mbox{SL}(3,{\bbbz})$ and three dimensional 
representions of $\mbox{Gal}(\overline{\bbbq}/\bbbq)$ (see \cite{GT1} and 
\cite{vierman}). 
This relation is the equality of local $L$-factors, for primes $p\leq 173$,
attached to the modular forms and to the Galois representation, see
Theorem \ref{main}.
The result gives some evidence  for general conjectures
of Langlands and Clozel \cite{C1}. 

The first three section follow closely
the notes from a seminar talk of the first
author at the s\'{e}minaire de th\'{e}orie des nombres de Paris
in January 1995. 
In the first section we briefly recall an instance of the relation between
elliptic modular forms and Galois representations. In the second section
we introduce the modular forms for $\mbox{GL}(3)$ and the Galois representations
are discussed in section three.

In section four we give some new examples of non-cusp forms for 
congruence subgroups of $\mbox{SL}(3,\bbbz)$ and we describe many of these in terms
of classical modular forms for congruence
subgroups of $\mbox{SL}(2,\bbbz)$. The last section deals with a Hodge theoretical 
aspect
of the algebraic varieties  (motives in fact) we used to define the Galois
representations.

It is a pleasure to thank Avner Ash, Kevin Buzzard, Bas Edixhoven and Jasper
Scholten, especially for their interest and help concerning Sect.~\ref{4.6}.

\section{Modular Forms: the $\mbox{GL}(2)$ Case }

Let $S_2(N)$ be the space of cusp forms of weight two for the 
congruence subgroup $\Gamma_0(N)\subset \mbox{SL}(2,\bbbz)$. 
Let $f=\sum a_ne^{2\pi i n}\in S_2(N)$ be a newform, thus $a_1=1$ and
$f$ is an eigenform for the Hecke algebra: $T_pf=a_pf$ for all prime 
numbers $p$
which do not divide $N$. 
For such a prime $p$ one defines the local $L$-factor of $f$ as
$$
L_p(f,s):=(1-a_pp^{-s}+p^{1-2s})^{-1},
$$
note that $L_p(f,s)$ is determined by the eigenvalue $a_p$.

In case all $a_p$ are in $\bbbz$, $f$ defines an
elliptic curve $E_f$, defined over $\bbbq$ ($E_f$ is a subvariety
of the Jacobian of the modular curve $X_0(N)$).
The Galois group $\mbox{Gal}(\overline{\bbbq}/\bbbq)$ acts on the $\ell^n$-torsion points 
of this curve which gives an $\ell$-adic representation:
$$
\rho_{f,\ell}:\mbox{Gal}(\overline{\bbbq}/\bbbq)\longrightarrow \mbox{GL}_2(\bbbq_\ell).
$$
The local $L$-factor of this representation for primes $p$ as above
does not depend on the choice of the prime $\ell\neq p$ and is defined by
$$
L_p(\rho_f,s):=\det(I-\rho_{f,\ell}(F_p)p^{-s})^{-1}
=(1- \mbox{trace} (\rho_{f,\ell}(F_p))p^{-s}+p^{1-2s})^{-1},
$$
with $F_p\in \mbox{Gal}(\overline{\bbbq}/\bbbq)$ a Frobenius element at $p$.

The Eichler-Shimura congruence relation asserts that
$$
a_p= \mbox{trace} (\rho_{f,\ell}(F_p))\qquad{\rm so}\quad
L_p(f,s)=L_p(\rho_f,s)
$$
(again with $p$ a prime not dividing $N\ell$).
Thus we have a method to associate to a newform $f$ a (compatible
system of $\ell$-adic) Galois representation(s) 
$\rho_{f,\ell}$ such that the $L$-factors agree.
This construction has been generalized to newforms of any weight (and 
arbitrary Hecke eigenvalues) by Deligne \cite{D}
using Galois representations on 
certain etale cohomology groups of certain $\ell$-adic sheaves
on the modular curve $X_0(N)$.

It is a pleasure to observe that recently Wiles proved a partial inverse
to the construction sketched above: he shows that
for a certain class of elliptic curves
defined over $\bbbq$ the corresponding Galois $L$-series are the $L$-series
of newforms. As is well known, this has been used to prove Fermat's Last
Theorem.

\section{Modular Forms for $\mbox{GL}(3)$}

\subsection{}
One can also define modular forms, a Hecke algebra and local $L$-factors 
for congruence subgroups of $\mbox{SL}(3,\bbbz)$, see below. 
However, the upper half plane 
$$
\bbbh=\{z\in\bbbc:\;\Im(z)>0\}\cong \mbox{SL}(2,\bbbr)/\mbox{SO}(2),
$$
which has a complex structure, is now replaced by $\mbox{SL}(3,\bbbr)/\mbox{SO}(3)$ (see 
\cite{AGG}),
a real variety of dimension 5 which, for dimension reasons(!), cannot have
a complex structure. 

In particular, one does not know how to associate
algebraic varieties to congruence subgroups of $\mbox{SL}(3,\bbbz)$ (in contrast
to the modular curves in the $\mbox{GL}(2)$-case). Therefore there
are no a priori given Galois representations on etale cohomology groups 
which could be related to modular forms for such congruence subgroups.

\subsection{}
In the case of $\mbox{SL}(2,\bbbz)$, the space of holomorphic modular
forms of weight two for a congruence subgroup $\Gamma$ is a subspace of the
cohomology group $H^1(\Gamma,\bbbc)$. This generalizes as follows.

\subsection{}
{}From now on we use the following definition:
$$
 \Gamma_0(N) = \Bigl\{ (a_{ij})\in \mbox{SL}(3,\bbbz) \mid
   a_{21} \equiv 0 \;{\rm mod}\; N \hbox{\rm\ and\ }
   a_{31} \equiv 0 \;{\rm mod}\; N \Bigr\}.
$$
The modular forms for $\Gamma_0(N)$ we consider are elements of
$H^3(\Gamma_0(N),\bbbc)$. 
To compute this vector space, we introduce a finite set:
$$
 \bbbp^2(\bbbz/N) = \left\{ (\bar x,\bar y,\bar z)\in(\bbbz/N)^3 \mid
   \bar x\bbbz/N + \bar y\bbbz/N + \bar z\bbbz/N = \bbbz/N \right\} \Big/
   (\bbbz/N)^\times .
$$
When the elements of this set are viewed as column vectors, there is a
natural left action of $\mbox{SL}(3,\bbbz)$ on $\bbbp^2(\bbbz/N)$.  This action
is transitive, and the stabilizer of $(\bar1\colon\bar0\colon\bar0)$
equals $\Gamma_0(N)$.  Therefore
$$
\mbox{SL}(3,\bbbz)/\Gamma_0(N) \cong \bbbp^2(\bbbz/N).
$$

This relation between $\Gamma_0(N)$ and $\bbbp^2(\bbbz/N)$ 
leads to a very concrete description
of the vector space $H^3(\Gamma_0(N),\bbbc)$. In fact, its dual  
$H_3(\Gamma_0(N),\bbbc)$ can be computed as follows:

\subsection{Theorem.}\label{space} (\cite{AGG}, Thm 3.2, Prop 3.12)

There is a canonical isomorphism between $H_3(\Gamma_0(N),\bbbc)$ and
the vector space of mappings $f:\bbbp^2(\bbbz/N)\to\bbbc$ that satisfy
\begin{enumerate}
\parskip2pt
\item   $f(\bar x\colon\bar y\colon\bar z)=
	     -f(-\bar y\colon\bar x\colon\bar z)$,
\item $f(\bar x\colon\bar y\colon\bar z)=
	     f(\bar z\colon\bar x\colon\bar y)$,
\item $f(\bar x\colon\bar y\colon\bar z)+
	     f(-\bar y\colon\bar x-\bar y\colon\bar z)+
	     f(\bar y-\bar x\colon-\bar x\colon\bar y)=0$.
\end{enumerate}

\subsection{} \label{hecke}
For any $\alpha\in \mbox{GL}(3,\bbbq)$ one has a ($\bbbc$-linear) Hecke
operator:
$$
T_\alpha:\,H^3(\Gamma_0(N),\bbbc)\longrightarrow H^3(\Gamma_0(N),\bbbc).
$$
The adjoint operator $T_\alpha^*$ on the dual space
$H_3(\Gamma_0(N),\bbbc)$ can be explicitly computed using modular symbols.

The Hecke algebra ${\cal T}$ is defined to be the subalgebra of
$\mathop{\rm End}
(H^3(\Gamma_0(N),\bbbc))$ generated by the $T_\alpha$'s with $\det(\alpha)$
relatively prime with $N$.
The Hecke algebra is a commutative algebra and  
we are interested in eigenforms $F\in H^3(\Gamma_0(N),\bbbc)$ for
the Hecke algebra:
$$
 TF = \lambda(T)F,\qquad{\rm with}\quad \lambda:{\cal T}\to\bbbc
 \qquad(\mbox{for all}\;T\in{\cal T}). 
$$

Of particular interest are the Hecke operators $E_p=T_{\alpha_p}$, 
which are for a prime $p$ not dividing $N$
defined using 
$
\alpha_p=\pmatrix{p&0&0\cr 0&1&0\cr 0&0&1\cr}\in \mbox{GL}(3,\bbbq)$.

Let
$a_p:=\lambda (E_p)$,
for a (given) character $\lambda$ of ${\cal T}$ and a prime $p$
not dividing $N$, then the local
$L$-factor of a Hecke eigenform $F\in H^3(\Gamma_0(N),\bbbc)$
(with the additional condition that $F$ is cuspidal) 
corresponding to $\lambda$ (so $E_pF=a_pF$) is
$$
L_p(F,s)=(1-a_pp^{-s}+\bar{a}_p p^{1-2s}-p^{3-3s})^{-1},
$$
where $\bar{a}_p$ is the complex conjugate of $a_p$. 
The field $K_F:=\bbbq(\ldots,a_p,\ldots)$ generated 
by the Hecke eigenvalues of an eigenform $F$ is known to be either totally
real or is a CM field (a degree 2, non-real extension of a totally real 
field).

\subsection{}\label{tabel}
In \cite{vierman}, a list of the $a_p$'s with $p\leq 173$ is given
 for several
eigenforms with $N\leq 245$. Here we list some 
$a_p$'s of three
particularly interesting eigenforms (these eigenforms are uniquely
determined by their level $N$ and the $a_p$'s listed). In case $p$ divides 
$N$ we write  $**$ for $a_p$.  In the three cases listed here $K_F=\bbbq(i)$ 
with $i^2=-1$. The 
complex conjugates of the $a_p$'s for a given $F$ are the Hecke eigenvalues
for another modular form $G$ of the same level.

{\small
$$
\begin{array}{|@{\hspace{3pt}}c||@{\hspace{3pt}}c|@{\hspace{3pt}}c|
@{\hspace{3pt}}c|@{\hspace{3pt}}c|@{\hspace{3pt}}c|@{\hspace{3pt}}c|
@{\hspace{3pt}}c|@{\hspace{3pt}}c|@{\hspace{3pt}}c|}
\hline p= & 2 & 3 & 5 & 7 & 11 & 13 & 17 & 101& 173 \\\hline \hline
N&\multicolumn{8}{c}{\mbox{eigenvalue $a_p$}}&\\ \hline
128 & ** & 1+2i & -1-4i & 1+4i & -7-10i & -1+4i & 7&-105-100i& -49-188i\\\hline
160 & ** & 1+2i & **    &1-2i & -3-12i& -5-8i & -5 &-33+64i& 99+104i \\\hline
205 & -1 & 1+2i & **    & 1+2i & -7-10i & 3-8i & -5 &115-40i& -153-288i \\\hline
\end{array}
$$
}

\section{Galois Representations}

\subsection{}
We are interested in relating Hecke eigenforms and Galois
representations. In particular, given a Hecke eigenform $F$ we would
like to find (a compatible system of) $\lambda$-adic Galois representations
$$
\rho_{F,\lambda}:\mbox{Gal}(\overline{\bbbq}/\bbbq)\longrightarrow \mbox{GL}(W_\lambda)
$$
having the same local $L$-factors as $F$.
Here $\lambda$ is a prime in a 
finite extension $K_\lambda$ of $\bbbq_\ell$ and $W_\lambda$ is a (finite
dimensional) $K_\lambda$ vector space. 
The local $L$-factors of $\rho_{F,\lambda}$ (again independent of $\lambda$)
being defined  as before (for unramified primes, conjecturally those not 
dividing
$N\ell$):
$$
L_p(\rho_F,s):=det(I-\rho_{F,\lambda}(F_p)p^{-s})^{-1}.
$$
In particular, we want $\dim W_\lambda=3$.

\subsection{} The case that $K_F$ is totally real is analyzed by Clozel 
\cite{C2}. We just recall that if in this case such a Galois 
representation $\rho_{F,\lambda}$ exists then $\rho_{F,\lambda}$
is selfdual in the following sense.

Consider the Tate-twisted dual Galois representation:
$$
\rho_{F,\lambda}^*:={}^t\rho_{F,\lambda}^{-1}(-2):
\mbox{Gal}(\overline{\bbbq}/\bbbq)\longrightarrow \mbox{GL}(W_\lambda),\qquad
{\rm so}\qquad
\rho_{F,\lambda}^*(F_p):=p^2\,{}^t\rho_{F,\lambda}^{-1}(F_p).
$$
Let $\alpha_i$, $i=1,\,2,\,3$ be the eigenvalues of
$\rho_{F,\ell}(F_p)$, then the eigenvalues of $\rho_{F,\ell}^*(F_p)$
are $\beta_i:=p^2/\alpha_i$. 
Since $\sum \alpha_i=a_p,\;\sum \alpha_i\alpha_j=pa_p$ (since now 
$\bar{a}_p=a_p$)
and $\prod \alpha_i=p^3$, the sets of eigenvalues $\{\alpha_i\}$ and 
$\{\beta_i\}$ coincide. 

Thus
$L_p(\rho_F,s)=L_p(\rho^*_F,s)$ for all $p$ not dividing $N$
and so the (semi-simplifications of the) Galois representations are 
the same. It implies also that a subgroup of finite index of the
image of $\mbox{Gal}(\overline{\bbbq}/\bbbq)$ is contained in a group 
$G\subset \mbox{GL}(W_\lambda)$ with $G\cong \mbox{PGL}(2,K_\lambda)$.
Examples of this are the $\mbox{Sym}^2$ of Galois representations in 
$\mbox{GL}(2,\bbbq_\ell)$. 

\subsection{}
We will be especially interested in the non-selfdual case. 
Since we found several examples of
Hecke eigenforms $F$ with $K_F=\bbbq(i)$ we will consider that case here.
To find corresponding Galois representations we use the fact that for
any algebraic variety $X$ defined over $\bbbq$, one has a Galois 
representation on the etale cohomology:
$$
\mbox{Gal}(\overline{\bbbq}/\bbbq)\longrightarrow 
\mbox{GL}(H^n_{et}(X_{\overline{\bbbq}},\bbbq_\ell)).
$$
The point is to find a suitable $X$ and (a subspace of) a suitable 
$H^n_{et}$.
In case $X$ is smooth, projective, and has good reduction mod $p$,
theorems of Grothendieck and Deligne imply that
the eigenvalue polynomial of $F_p$ acting on 
$H^n_{et}(X_{\overline{\bbbq}},\bbbq_\ell)$
has coefficients in $\bbbz$, is independent
of $\ell$ and the eigenvalues of $F_p$ have absolute value $p^{n/2}$.

 The desired
equality $L_p(F,s)=L_p(\rho_F,s)$
for the eigenforms $F$ from (\ref{tabel})
(and one expects the same more generally
for certain cusp forms, `Ramanujan conjecture'),
implies that the absolute value of the eigenvalues of $\rho_F(F_p)$
must be $p$. Therefore we will consider $H^2_{et}$ and take $\dim X>1$
since $\dim H^2_{et}=1$ for curves. 

A well-known theorem 
implies that $H^2_{et}(X_{\overline{\bbbq}},\bbbq_\ell)\hookrightarrow 
H^2_{et}(S_{\overline{\bbbq}},\bbbq_\ell)$
where $S$ is a suitable surface contained in $X$. Thus we restrict 
ourselves to considering $H^2_{et}(S_{\overline{\bbbq}},\bbbq_\ell)$ for a
surface $S$.

The Galois representation on this $\bbbq_\ell$-vector space is reducible
in general, a decomposition is:
$$
H^2_{et}(S_{\overline{\bbbq}},\bbbq_\ell)=
T_\ell\oplus \,\mbox{NS}(S_{\overline{\bbbq}})\otimes_\bbbz\bbbq_\ell
$$
where $\mbox{NS}(S_{\overline{\bbbq}})$ is the N\`eron-Severi group of the surface
$S$ over $\overline{\bbbq}$ (the Galois group  permutes the classes of 
divisors modulo a Tate twist) and $T_\ell$ is the orthogonal complement of
$\mbox{NS}(S_{\overline{\bbbq}})$ with respect to the intersection form.
The intersection form is the cup product $H^2_{et}\times H^2_{et}\rightarrow
H^4_{et}\cong \bbbq_\ell$. The eigenvalues of Frobenius on 
$\mbox{NS}(S_{\overline{\bbbq}})\otimes\bbbq$ are roots of unity multiplied by $p$,
so $\rho_{F,\lambda}$, if it exists, should be a representation on 
a subspace of 
$T_\ell\otimes_{\bbbq_\ell}K_\lambda$.

In case $T_\ell$ has dimension 3, the Galois representation on it will be 
selfdual (due to the intersection form).
To find a 3 dimensional
Galois representations with $trace{F_p}\in\bbbz[i]$ as desired
we assume that the surface has
an automorphism, defined over $\bbbq$:
$$
\phi:S\longrightarrow S,\qquad{\rm with}\quad \phi^4=id_S.
$$
Thus $\phi^*:H^2_{et}\rightarrow H^2_{et}$ will commute with the Galois 
representation.

Assume moreover that $\dim T_\ell=6$ and 
$\phi^*:T_\ell\rightarrow T_\ell$ has two 3-dimensional
eigenspaces $W_\lambda,\;
W'_\lambda$ (with eigenvalue $\pm i$):
$$
T_\lambda:=T_\ell\otimes_{\bbbq_\ell}K_\lambda=W_\lambda\oplus
W'_\lambda
$$
with $K_\lambda$ an extension of $\bbbq_\ell$ containing $i$. Then we have
a 3-dimensional Galois representation $\sigma'$ on $W_\lambda$.
The determinant of $\sigma'(F_p)$
is in general not equal to $p^3$ but is $\chi(p)p^3$ for a Dirichlet
character $\chi$. Twisting $\sigma'$ by this character
we get a Galois representation 
$$
\sigma_{S,\lambda}:\mbox{Gal}(\overline{\bbbq}/\bbbq)\longrightarrow \mbox{GL}(W_\lambda).
$$
whose $L$-factors $L_p(\sigma_S,s)$ are similar to the $L_p(F,s)$ for the 
eigenforms in the example above.

Note that the intersection form $(\cdot\,,\cdot)$
restricted to $W_\lambda$ is trivial (because it
is invariant under pull-back by $\phi^*$ and extends $K_\lambda$-linearly:
$(w_1,w_2)=(\phi^*w_1,\phi^*w_2)$ $=(iw_1,iw_2)=i^2(w_1,w_2)$ $=-(w_1,w_2)$ with
$w_1,\,w_2\in W_\lambda$). Thus there is no obvious reason for $\sigma_S$ to be
selfdual.

\subsection{}
Now one has to search for such surfaces. The main problem is that in 
general $\dim H^2_{et}$ will be large but $\mbox{rank} \mbox{NS}$ will be small.
Thus it is  not so easy to get $\dim T_\ell=6$,
see however \cite{GT1} and \cite{GT2} for various examples.

The most interesting example is given by the one parameter family of surfaces
$S_a$ which are the smooth, minimal, projective model of the 
singular, affine surface defined in $x,\,y,\,t$-space by
$$
t^2=xy(x^2-1)(y^2-1)(x^2-y^2+axy),\qquad{\rm and}\quad
 (x,y,t)\longmapsto (y,-x,t)
$$
defines the automorphism $\phi$.
In \cite{GT1}, 3.7-3.9, we explain how to determine eigenvalue polynomials
of $\sigma_{S,\lambda}(F_p)$, and thus the $L$-factors,
basically using the Lefschetz trace formula and
counting points on $S$ over finite fields.
The main result is:

\subsection{Theorem.}\label{main}(\cite{GT1}, 3.11; \cite{vierman}, 3.9)
The local $L$-factors of the modular forms for $N=128,\,160,\,205$
in $\S$\ref{tabel} are the
same as the local $L$-factors of the Galois representations 
$\sigma_{S_a,\lambda}$,
with $a=2,\,1,\,\mbox{$\frac{1}{16}$}$ respectively,
for all odd primes $p\leq 173$ not dividing $N$.

\subsection{}
In \cite{GT2} we gave
another construction of surfaces $S$ which define 3 dimensional Galois
representations.
 These surfaces are degree $4$ cyclic base changes of elliptic surfaces
${\cal E}\rightarrow \bbbp^1$.
By taking the orthogonal complement to a large algebraic part in
$H^2_{et}$ together with all cohomology coming from the intermediate degree $2$ 
base
change, one obtains a representation
space, similar to $T_\ell$, for $\mbox{Gal}(\overline{\bbbq}/\bbbq)$. 
Taking an eigenspace $W_\lambda$
of  the action of the automorphism of order $4$ defining the cyclic base
change finally gives 3 dimensional Galois representations. 

Our main (technical) result is a formula for the traces of Frobenius elements
on this space in terms of the number of points on
${\cal E}$ and $S$ over a finite
field (\cite{GT2},Theorem 3.4). 
This formula allows us to compute the characteristic
polynomial
of Frobenius in many cases. 

We use this result to prove that certain
examples yield selfdual representations, while others do not.
For some of the selfdual cases we can actually
exhibit 2-dimensional Galois representations (defined by elliptic curves)
whose symmetric square seems 
to
coincide with the 3-dimensional Galois representation. 

We did not find new examples of non-selfdual Galois representations with
the same local $L$-factors as modular forms, probably because the
conductor of these Galois representations is too large.
We would like to point out that there does not seem to be an explicit way
to determine the conductor of the Galois representation $\sigma_S$
in terms of the geometry of $S$ (a surface over $\mbox{Spec}(\bbbq)$).

\section{Non-Cusp Forms and Galois Representions}

\subsection{}
In this section we give an example of the decomposition in Hecke eigenspaces
of a cohomology group $H^3(\Gamma_0(N),\bbbc)$.
We will take $N=245$. This
example is also mentioned in \cite{vierman}, $\S$3.5 where
it is shown that a certain 8 dimensional Hecke invariant
subspace of $H^3(\Gamma_0(245),\bbbc)$ contains no cusp forms.
Here we extend this by interpreting most of the $83$ dimensional
space $H^3(\Gamma_0(245),\bbbc)$ in terms of so-called
Eisenstein liftings of classical elliptic cusp forms and
of Eisenstein series. 

As before, if $F\in H^3(\Gamma_0(245),\bbbc)$ is an eigenform
for all Hecke operators, we denote by $K_F$ the field
generated by all eigenvalues of the Hecke operators on $F$.
As a first step towards the decomposition we have the following
Proposition.

\begin{proposition}\label{dims}
The cohomology group $H^3(\Gamma_0(245),\bbbc)$ decomposes as
$$
H^3(\Gamma_0(245),\bbbc)=V_1\oplus V_2\oplus V_3\oplus V_4\oplus V_5
$$
(as a module over the Hecke algebra),
with
\begin{itemize}
\item $\dim V_1=25$ and $V_1$ is generated by eigenforms $F$
with $K_F=\bbbq$;
\item $\dim V_2=16$ and $V_2$ is generated by eigenforms $F$
with $K_F=\bbbq(\sqrt{2})$;
\item $\dim V_3=16$ and $V_3$ is generated by eigenforms $F$
with $K_F=\bbbq(\sqrt{17})$;
\item $\dim V_4=8$ and $V_4$ is generated by eigenforms $F$
with $K_F=\bbbq(\sqrt{2},\sqrt{-3})$;
\item $\dim V_5=18$ and $V_5$ is generated by eigenforms $F$
with $K_F=\bbbq(\sqrt{-3})$.
\end{itemize}
None of the spaces $V_1,\ldots,V_5$ contains a non-zero cuspform;
in fact, these spaces are generated by Eisenstein liftings or
(in the case of $V_4$ and $V_5$) twists of such by cubic Dirichlet characters.

\end{proposition}

\subsection{}
With notations as given in \cite{vierman} $\S$3.5, one has
$V_4=V_a\oplus V_b$, hence this case of the above
proposition is already described in {\sl loc. sit.}

We briefly recall the two types of Eisenstein liftings
of classical modular forms here. Let $f$ be a normalized
elliptic cuspform of level $N$ and weight 2, which is an 
eigenform for
the Hecke operators $T_n$ with $(n,N)=1$. Also, we allow
$f$ to be the normalized Eisenstein series of weight 2:
$f=-B_2/4+\sum_{n=1}^{\infty}\sigma_1(n)q^n$; so $a_p=p+1$,
compare e.g.
\cite{Kob} for notations. The Fourier coefficients
in the $q$-expension $f=q+a_2q^2+a_3q^3+\cdots$ define a
Dirichlet series $L(f,s)=\sum_n a_n n^{-s}$. This series has an
Euler product expansion with Euler factors 
$(1-a_pp^{-s}+p^{1-2s})^{-1}$ for primes $p$ which
do not divide $N$ 
(in case $f$ is the Eisenstein series, these factors are
$(1-p^{-s})^{-1}(1-p^{1-s})^{-1}$).

Given $f$, one constructs two eigenclasses $F_1,F_2\in
H^3(\Gamma_0(N),\bbbc)$. The $F_1$ has eigenvalue $pa_p+1$
for the $p$th Hecke operator $E_p$, and $F_2$ eigenvalue
$a_p+p^2$. 

On the Galois side of the Langlands correspondence, it is
relatively easy to describe these liftings. Namely,
if $f$ corresponds to a 2 dimensional $\lambda$-adic 
representation space $V$ for $\mbox{Gal}(\overline{\bbbq}/\bbbq)$,
then $F_1$ corresponds to 
$V(-1)\oplus \bbbq_\lambda(0)$ and
$F_2$ to 
$V\oplus \bbbq_\lambda(-2)$ where $\bbbq_\lambda(n)$ is the 1 dimensional
$\lambda$-adic representation space on which the Galois group acts 
by the $-n$-th power of the cyclotomic character (thus $F_p$ acts as $p^{-n}$).
In case $f$ is the Eisenstein series, we have 
$V=\bbbq_\lambda(0)\oplus\bbbq_\lambda(-1)$ and the two lifted
representations coincide (both are $\bbbq_\lambda(0)\oplus
\bbbq_\lambda(-1)\oplus\bbbq_\lambda(-2)$).

\subsection{}
There exists a unique normalized cuspform of weight 2
and level $35$ which has $\bbbq$-rational Fourier coefficients.
This form yields $2$ eigenclasses in $H^3(\Gamma_0(35),\bbbc)$;
from the theory of oldforms \cite{Ree}, each of these
appears three times at level $35\cdot 7=245$.

Similarly, the modular form corresponding to the CM
elliptic curve of conductor $49$ gives rise to six
oldforms which are Eisenstein liftings.

Starting from the Eisenstein series, one finds $7$
forms at level $245$ all with eigenvalues $1+p+p^2$.

Finally, from tables of Cremona (as well as from
unpublished tables of Cohen, Skoruppa and Zagier) it
follows that there exist $3$ (elliptic) newforms of level $245$
which are Hecke eigenforms with rational eigenvalues.
Each of them gives us two Eisenstein liftings.

Adding up, we now have $6+6+7+6=25$ eigenclasses of level $245$
with rational eigenvalues. Our calculations made for
the tables in \cite{vierman} revealed that, e.g.,
the Hecke operator $E_2$ has precisely $25$ rational
eigenvalues (counted with multiplicity). Hence the
conclusion is, that the space $V_1$ given in Proposition~\ref{dims}
indeed has $\dim V_1=25$, and it is generated by
Eisenstein liftings as claimed.

\subsection{}
The cases $V_2,V_3$ are completely
analogous. For $V_2$, we note that there exist newforms of
weight $2$ and level $245$ with $q$-expansion
$q+\sqrt{2}q^2+(1+\sqrt{2})q^3+\ldots$ and
$q+(1+\sqrt{2})q^2+(1-\sqrt{2})q^3+\ldots$ respectively.
These together with their Galois conjugate forms
and their twists by the quadratic Dirichlet character modulo $7$
give us $8$ newforms of level $245$. Each of them yields two
Eisenstein liftings, and this precisely describes the space $V_2$
of dimension $16$.

Similarly, there are exactly two (conjugate) newforms of level $35$
with Fourier coefficients generating $\bbbq(\sqrt{17})$. They
provide $2\cdot 2=4$ Eisenstein liftings of level $35$, and hence
$3\cdot 4=12$ oldforms of level $245$. Twisting the newforms
by the quadratic character modulo $7$ yields newforms of
level $245$, and from these we find another $4$ Eisenstein liftings.
In this way, $V_3$ is generated.

\subsection{}\label{4.6}
Having described $V_1,\ldots,V_4$ (the latter space was already
treated in \cite{vierman}), and observing from Table~3.3
that $\dim H^3(\Gamma_0(245),\bbbc)=83$, we conclude we still
have to describe a Hecke-invariant space of dimension
$83-(25+16+16+8)=18$. To this end, we mention that at level
$49=245/5$, our programs found a $6$ dimensional Hecke invariant
subspace on which the operator $E_2$ acts with $6$ (pairwise
conjugate, pairwise different) eigenvalues in $\bbbq(\sqrt{-3})$.
Hence this space yields eigenforms with $K_F=\bbbq(\sqrt{-3})$.
Moreover, it lifts to a Hecke invariant subspace of dimension
$3\cdot 6=18$ at level $245$, which therefore exactly equals
the summand $V_5$ of $H^3$ we did not describe yet.

As an example, the eigenvalues of
the operator $E_3$ on $V_5$ are $a_3, \overline{a_3}, a_3\omega$,
$\overline{a_3\omega}, a_3\overline{\omega}$ and
$\overline{a_3}\omega$ where $\omega^2+\omega+1=0$ and
$a_3=-5-3\sqrt{-3}$. This situation is explained as
follows. The Euler factor that corresponds to a Hecke
eigenclass is obtained using the polynomial
$X^3-a_pX^2+pb_pX-p^3$, where $a_p$ is the eigenvalue
of the operator $E_p$. The number $b_p$ similarly corresponds
to the operator $D_p=T_{\beta_p}$, defined using $\beta_p
:=\pmatrix{p&0&0\cr 0&p&0\cr 0&0&1\cr}\in \mbox{GL}(3,\bbbq)$.
If the eigenclass is cuspidal, then $b_p$ is 
the complex conjugate of $a_p$. This in fact follows from
the fact that the associated automorphic representation
is unitary in that case. In our situation however,
a computation shows that $b_3=a_3\neq \overline{a_3}$.
Hence the representation cannot be unitary and therefore
the eigenclasses here are not cuspidal.

Based on calculations for primes $\leq 131$, the Hecke
eigenvalues seem to be as follows. For $p\neq 5,7$ 
we have $b_p=a_p=
\chi(p)( \psi(p) + p + \psi^2(p)p^2)$
with $\chi,\psi$  Dirichlet characters modulo $7$ of order
dividing $3$. 
This corresponds to the sum of 1-dimensional
Galois representations
$$(\chi\psi\otimes\bbbq_\lambda(0))\oplus(\chi\otimes
\bbbq_\lambda(-1))\oplus(\chi\psi^2\otimes\bbbq_\lambda(-2)).$$

\section{Variations of Hodge Structures of Weight Two}

\subsection{}
In all our constructions for Galois representations
we consider a subspace $T_\ell\subset 
H^2_{et}(S_{\overline{\bbbq}},\bbbq_\ell)$. 
This subspace is defined using algebraic cycles, 
thus there exists also a Betti realization $T_\bbbz\subset 
H^2(S(\bbbc),\bbbz)$ (of the motive $T$) which is a polarized Hodge structure
of weight two. We recall the relevant definitions and the main results
of Griffiths and Carlson on the moduli of the $T_\bbbz$'s. 

The main point is the essential difference with the weight one case
(which is essentially the theory of abelian varieties). In the weight one
case, one has a universal family of abelian varieties over suitable
quotients of the Siegel space. In the weight two (and higher) case, the
analogy of the Siegel space is a certain (subset of a) period domain,
but in general (and in particular this is the case with the $T_\bbbz$ under
consideration), the (polarized) Hodge structures obtained from algebraic varieties
do not fill up the period space. In fact we will see that the 
Hodge structures like $T_\bbbz$ are parametrized by a 4-dimensional space,
but those that come from geometry have at most a 2-dimensional deformation
space (and imposing an automorphism of order 4 as we do 
implies a 1-dimensional deformation space). 

It is not clear whether these
period spaces (or the subvarieties parametrizing `geometrical' Hodge 
structures) have good arithmetical properties like Shimura varieties.

\subsection{}
Recall that a $\bbbz$-Hodge structure $V$ of weight $n$ is a free $\bbbz$-module of 
finite rank together with  decomposition:
$$
V_\bbbc:=V\otimes_\bbbz\bbbc=\oplus_{p+q=n}V^{p,q},\qquad
{\rm with}\quad
\overline{V^{p,q}}=V^{q,p},
$$
where the $V^{p,q}$ are complex vector spaces and the bar indicates complex
conjugation (given by $\overline{v\otimes z}=v\otimes\overline{z}$).

A rational Hodge structure $V_\bbbq$ is a finite dimensional $\bbbq$-vector space 
with a similar decomposition of $V_\bbbc:=V_\bbbq\otimes_\bbbq\bbbc$. Thus a 
$\bbbz$-Hodge structure $V$ determines a rational Hodge structure on 
$V_\bbbq:=V\otimes_\bbbz\bbbq$.

A (rational) Hodge structure $V_\bbbq$ determines an $\bbbr$-linear map,
the Weil operator:
$$
J:V_\bbbr:=V_\bbbq\otimes_\bbbq\bbbr\longrightarrow V_\bbbr\qquad
{\rm with}\quad
J_\bbbc v_{p,q}=i^{p-q}v_{p,q}
$$
for all $v_{p,q}\in V^{p,q}$ and $J_\bbbc$ is the $\bbbc$-linear extension of 
$J$. 
One has $J^2=(-1)^n$ since $i^{2p-2q}=(-1)^{p-q}=(-1)^{p+q}$.
Thus $J$ determines a complex structure on 
$V_\bbbr$ in case $V$ has odd weight.

A polarization on a rational Hodge structure $V_\bbbq$ of weight $n$ is a 
bilinear 
map
$$
\Psi:V_\bbbq\times V_\bbbq\longrightarrow \bbbq,\qquad
\Psi_\bbbc(v_{p,q},v_{r,s})=0\quad{\rm unless}\quad p+r=q+s=n
$$
(intrinsically: $\Psi:V_\bbbq\otimes  V_\bbbq\rightarrow \bbbq(-n)$ is a morphism of 
Hodge  
structures)
which satisfies the Riemann relations, that is, for all $v,\,w\in V_\bbbr$: 
$$
\Psi(v,Jw)=\Psi(w,Jv),\qquad
\Psi(v,Jv)>0\quad ({\rm if}\; v\neq 0)
$$
thus $\Psi$ defines an inner product $\Psi(-,J-)$ on $V_\bbbr$.

One easily verifies, using the first property, that 
$\Psi(Jv,Jw)=\Psi(v,w)$, since also 
$\Psi(Jv,Jw)=\Psi(w,J^2v)=(-1)^n\Psi(w,v)$, a polarization is symmetric
if $n$ is even and antisymmetric if $n$ is odd.

\subsection{}
For a smooth complex projective variety $X$ the cohomology groups 
$H^n(X,\bbbq)$ are polarized rational Hodge structures of weight $n$.
One writes $H^{p,q}(X):=H^n(X,\bbbc)^{p,q}$. 
In case $X$ is a surface, the cup product on $H^2(X,\bbbq)$ 
 (note that $H^4(X,\bbbq)=\bbbq$) gives ($-1$ times)
 a polarization on the primitive cohomology
 $H^2_{prim}$. In particular it induces a polarization $\Psi$
 on the sub-Hodge structure 
 $T_\bbbq=\mbox{NS}^{\perp}$
 of $H^2(S(\bbbc),\bbbq)$ which we consider.

\subsection{}\label{polT} 
Let $T_\bbbz$ be a Hodge structure of weight 2 and rank 6 with
$$
T_{\bbbc}=T^{2,0}\oplus T^{1,1}\oplus T^{0,2},\qquad \dim T^{p,q}=2
$$
for all $p,\,q$. 
Then one easily verifies that:
$$
T_{\bbbr}=W_1\oplus W_2\qquad{\rm with}\quad
\left\{\begin{array}{rcl}
W_1&:=&T_\bbbr\cap T^{1,1}\\
W_2&:=&T_\bbbr\cap (T^{2,0}\oplus T^{0,2})
\end{array}\right.
$$

For $v\in W_1\subset T^{1,1}$ we have $Jv=v$ and thus 
$\Psi(v,v)=\Psi(v,Jv)>0$, so $\Psi$ is positive definite on $W_1$.
Hence we can choose
an $\bbbr$ basis $f_1,\,f_2$ of $W_1$
which is orthonormal w.r.t.\ $\Psi$ and which is a $\bbbc$-basis of 
$T^{1,1}=W_1\otimes_\bbbr\bbbc$.

For $v\in W_2$ we have $v=v_{2,0}+v_{0,2}$ thus $Jv=-v$ and so $\Psi$
is negative definite on $W_2$. Let $v_1:=e_1+\bar{e}_1,\quad 
v_2:=e_2+\bar{e}_2$ be an orthonormal basis for $(-1/2)\Psi$ on $W_2$ with
$e_1,\,e_2\in V^{2,0}$.
Then $e_1,\,e_2$ is a $\bbbc$-basis of $T^{2,0}$ (and thus 
$\bar{e}_1,\,\bar{e}_2$ is a $\bbbc$-basis of $T^{0,2}$).
Note $-2=\Psi(e_1+\bar{e}_1,e_1+\bar{e}_1)=
\Psi(e_1,\bar{e}_1)+\Psi(\bar{e}_1,e_1)=2\Psi(e_1,\bar{e}_1)$
(since $\Psi$ is symmetric).
In this way one finds
$\Psi(e_k,\bar{e}_l)=-\delta_{kl}$ (Kronecker's delta) thus $\Psi_\bbbc$
is given by the matrix $Q$ on the basis $e_1,\,e_2,\,f_1,\,f_2,\,
\bar{e}_1,\,\bar{e}_2$ of $T_\bbbc$:
$$
Q=\left(\begin{array}{ccc}0&0&-I\\ 0&I&0\\-I &0&0 \end{array}\right).
$$

\subsection{} We consider first order deformations of the polarized
Hodge structure $T_\bbbz$ as in $\S$\ref{polT}.
Thus we fix the $\bbbz$-module and the 
bilinear map $\Psi$ and consider deformations of the Hodge structure
induced by deformations of an algebraic variety $X$ with $T_\bbbz\subset H^2(X,\bbbz)$, that
is, of the direct sum decomposition $T_\bbbc=\oplus T^{p,q}$. 

The first order deformations of a smooth complex projective algebraic variety
$X$ are parametrized by $H^1(X,\Theta_X)$
with $\Theta_X$ the tangent bundle of $X$
(Kodaira-Spencer theory). The isomorphisms $H^{p,q}(X)=H^q(X,\Omega^p)$
and the contraction map 
$\Theta_X\otimes_{{\cal O}_X}\Omega_X^p\rightarrow \Omega_X^{p-1}$ give a 
cup product map:
$$
H^1(X,\Theta_X)\otimes H^{p,q}(X)\longrightarrow H^{p-1,q+1}(X).
$$
Thus, for any $n$, we obtain a map, called the infinitesimal period map:
$$
\delta:H^1(X,\Theta_X)\longrightarrow 
\oplus_{p+q=n} {\rm Hom}(H^{p,q}(X),H^{p-1,q+1}(X)).
$$
Griffiths proved that for $\theta\in H^1(X,\Theta_X)$, the deformation
of the Hodge structure induced by the deformation of $X$ in the direction of $\theta$ is essentially given by $\delta(\theta)$.


The subspace $\Im(\delta)$
of $\oplus_{p+q=n}{\rm 
Hom}(H^{p,q}(X),H^{p-1,q+1}(X))$
satisfies (at least) two conditions. The first comes from
the polarization (see $\S$\ref{polcon}), the second is an integrability
condition found by Griffiths which is non-trivial only if the weight
of the Hodge structure is greater than one (see $\S$\ref{grifcon}).

We will now spell out the restriction of these conditions to the 
sub Hodge structure
$\oplus_{p+q=n}{\rm Hom}(T^{p,q},T^{p-1,q+1})$.

\subsection{}\label{polcon}
The condition that $\psi\in\oplus_{p+q=n}{\rm Hom}(T^{p,q},T^{p-1,q+1})
\subset \mbox{End}(T_\bbbc)$
preserves the polarization on $T$, is
 that $\Psi((I+t\psi)v,(I+t\psi)w)=\Psi(v,w)$ when 
$t^2=0$:
$$
\Psi_\bbbc(\psi(v),w)+\Psi_\bbbc(v,\psi(w))=0\qquad
\forall\;x,\,y\in T_\bbbc
$$
This condition implies that if $\psi$ preserves $\Psi$, then
it is determined by $\psi_{2}$ where
$$
\psi=(\psi_{2},\psi_{1})\in {\rm Hom}(T^{2,0},T^{1,1})\oplus
{\rm Hom}(T^{1,1},T^{0,2}).
$$
In fact,
 for all $v\in T^{2,0}$ and $w\in T^{1,1}$ we now have:
$
\Psi_\bbbc(v,\psi_1(w))=-\Psi_\bbbc(\psi_2(v),w)$.
Since $\Psi_\bbbc$ identifies $(T^{0,2})^{dual}$ with $T^{2,0}$, this 
equality thus defines $\phi_1(w)$ in terms of $\phi_2$. 
 
\subsection{}\label{polpres}
 With respect to the basis of $T_\bbbc$ considered in 
\ref{polT}, $\psi\in {\rm Hom}(T^{2,0},T^{1,1})\oplus
{\rm Hom}(T^{1,1},T^{0,2})\subset \mbox{End}(T_\bbbc)$ is given by a matrix
$N$ and the condition on $\psi$ becomes ${}^tNQ+QN=0$ so:
$$
N=\left(\begin{array}{ccc}0&0&0\\ A&0&0\\0&B&0\end{array}\right)
\quad{\rm and}\quad 
B={}^tA
$$
where the matrix $A$ (defining $\phi_2:T^{2,0}\rightarrow T^{1,1}$) can be 
chosen arbitrarily. This gives an isomorphism between the space $M_2(\bbbc)$
of $2\times 2$ complex matrices and polarization preserving deformations
$\psi$:
$$
M_2(\bbbc)\stackrel{\cong}{\longrightarrow }({\rm Hom}(T^{2,0},T^{1,1})\oplus
{\rm Hom}(T^{1,1},T^{0,2}))_\Psi,
\quad
A\longmapsto N(A):=
\left(\begin{array}{ccc}0&0&0\\ A&0&0\\0&{}^tA&0\end{array}\right).
$$
Thus we have a four dimensional deformation space. In case of Hodge structures 
of weight one, preserving the polarization is the only infinitesimal
condition. Here, in the weight two case, there is however another condition.

\subsection{}  \label{grifcon}
An important restriction, discovered by Griffiths, on the image of $\delta$ is:
$$
[{\rm Im}\delta,\,{\rm Im}\delta]=0\qquad{\rm i.e.}\quad
\delta(\alpha)\circ\delta(\beta)=\delta(\beta)\circ\delta(\alpha),
$$
for all $\alpha,\,\beta\in H^1(X,\Theta_X)$, 
so $\Im(\delta)$ is an abelian subspace of 
$\mbox{End}(T_\bbbc)$. For Hodge structures of weight 
$n\geq 2$ this imposes non-trivial conditions on the (dimension of) the 
image 
of $\delta$. We consider again our example (cf.\ \cite{Car}).

\subsection{}
We already determined the polarization preserving deformations in
$\S$\ref{polpres}. Using the same notation we find that Griffiths' condition 
is:
$$
N(A)N(B)=N(B)N(A)\qquad{\rm thus}\quad {}^tAB={}^tBA.
$$
This condition can be rephrased as saying that ${}^tAB$ must be symmetric.

Thus the image of $\delta$ is at most two dimensional and if
it is two dimensional with basis $N(A),\;N(B)$ then $A$ and $B$
span a maximal isotropic subspace
of the symplectic form:
$$
E:M_2(\bbbc)\times M_2(\bbbc)\longrightarrow \bbbc,
\quad E(A,B):=a_{11}b_{12}-a_{12}b_{11}+a_{21}b_{22}-a_{22}b_{21}=0.
$$
We recall that we also have an automorphism $\phi^*:T\rightarrow T$, 
preserving this automorphism gives another non-trivial condition on the
deformations. Thus the one parameter in our surfaces $S_a$ (and in the
other examples from \cite{GT2}) is the maximal possible.

\end{document}